\documentclass[10pt]{article} 
\usepackage[left=1in,top=1in,right=1in,bottom=1in]{geometry}
\usepackage{amsfonts,amssymb,amsmath,latexsym,indentfirst,} 
 
\newtheorem{rema}{Remark}[section] 
\newtheorem{defi}{Definition}[section] 
\newtheorem{lemm}{Lemma}[section] 
\newtheorem{theo}{Theorem}[section]

\newcommand{\Z}[1][]{\ensuremath{{\mathbb{Z}^{#1}} }} 

\newcommand{\C}[1][]{\ensuremath{{\mathbb{C}^{#1}} }} 
\newcommand{\R}{\ensuremath{{\mathbb{R}} }} 
\newcommand{\T}{\ensuremath{{\mathbb{T}} }}

\newcommand{\peq}{\hspace*{0.10in}}
\newcommand{\ppeq}{\hspace*{0.05in}}

\newcommand{\fim}{\rightline{$\blacksquare$}}

\author{Luiz Gustavo Farah\footnote{Partially supported by FAPESP-Brazil.} and M\'arcia Scialom.\\
\\
Department of Mathematics - IMECC/UNICAMP\\
C.P. 6065, Campinas, SP, BRAZIL\\
CEP 13083-970}
\title{On the periodic ``good'' Boussinesq equation \footnote{Mathematical subject classification: 35B30, 35Q55, 35Q72.}} 
\date{}

\begin{document} 
\maketitle 

\begin{abstract} 
We study the well-posedness of the initial-value problem for the periodic nonlinear ``good'' Boussinesq equation. We prove that this equation is local well-posed for initial data in Sobolev spaces \textit{$H^s(\T)$} for $s>-1/4$, the same range of the real case obtained in Farah \cite{LG4}.
\end{abstract} 

\section{Introduction} 
In this work we consider periodic boundary value problem (PBVP) for the Boussinesq-type equation
\begin{eqnarray}\label{NLB}
\left\{ 
\begin{array}{l}
u_{tt}-u_{xx}+u_{xxxx}+(f(u))_{xx}=0, \peq  x\in \T, t>0,\\
u(0,x)=u_0(x); \ppeq u_t(0,x)=u_1(x).
\end{array} \right. 
\end{eqnarray}

Equations of this type, in the continuous case, were originally derived by Boussinesq \cite{BOU} in his study of nonlinear, dispersive wave propagation.
We should remark that it was the first equation proposed in the literature to describe this kind of physical phenomena. 


Our principal aim here is to study the local well-posedness (LWP) for the PBVP associated to the ``good'' Boussinesq equation, that is, $f(u)=u^2$ in equation (\ref{NLB}), under low regularity of the data. Natural spaces to measure this regularity are the classical Sobolev spaces $H^s(\T)$, $s\in \R$ equipped with the norm
\begin{equation*}
\|f\|_{H^{s}(\T)}=\|\langle n\rangle^s \widehat{f}\|_{L^{2}(\T)}
\end{equation*}
where $\langle a\rangle\equiv 1+|a|$.

Concerning the LWP question for the initial value problem (IVP), several results have been obtained for the equation (\ref{NLB}).  Using Kato's abstract theory for quasilinear evolution equation, Bona and Sachs \cite{BS} showed LWP for $f\in C^{\infty}$ and initial data $u_0 \in H^{s+2}(\R)$, $u_1 \in H^{s+1}(\R)$ with $s > \frac{1}{2}$. Tsutsumi and Matahashi \cite{TM} established similar result when $f(u)=|u|^{p-1}u$, $p>1$ and $u_0 \in H^{1}(\R)$, $u_1 = \chi_{xx}$ with $\chi \in H^{1}(\R)$. Later, Linares \cite{FL} proved that (\ref{NLB}) is locally well-posedness in the case $f(u)=|u|^{p-1}u$, $1<p<5$ and $u_0 \in L^{2}(\R)$, $u_1 = h_{x}$ with $h \in H^{-1}(\R)$. The main tool used in his argument was the use of Strichartz estimates satisfied by solutions of the linear problem. We should remark that all these results also hold for the ``good'' Boussinesq equation. These results were improved by Farah \cite{LG4} who proved that the ``good'' Boussinesq equation is locally well-posedness for initial data and $(u_0,u_1) \in H^{s}(\R)\times H^{s-1}(\R)$, with $s>-1/4$

In the periodic setting, equation (\ref{NLB}) was studied by Fang and Grillakis \cite{FG}. Using the Fourier restriction norm approach introduced by Bourgain \cite{B} in his study of the nonlinear Schr\"odinger equation (NLS)
\begin{equation*}
iu_{t}+u_{xx}+u|u|^{p-2}=0, \textrm{ with } p\geq 3,
\end{equation*}
they proved LWP for (\ref{NLB}) assuming $u_0 \in H^{s}$, $u_1 \in H^{s-2}$, with $0\leq s\leq 1$ and  $|f(u)|\leq c|u|^{p}$, with $1<p<\frac{3-2s}{1-2s}$ if $0\leq s<\frac{1}{2}$ and $1<p<\infty$ if $\frac{1}{2}\leq s\leq 1$. Moreover, if $u_0 \in H^{1}$, $u_1 \in H^{-1}$ and  $f(u)= \lambda|u|^{q-1}u-|u|^{p-1}u$, with $1<q<p$ and $\lambda \in \R$ then the solution is global. In particular, for the ``good" Boussinesq equation, the lower Sobolev index where they obtain LWP is $s=0$.

In this paper, we also consider the periodic setting and improve this last result by establishing the LWP with $s> -1/4$ for the ``good'' Boussinesq equation. To this end we follow the argument in \cite{LG4}, where Farah use the method developed by Kenig, Ponce and Vega \cite{KPV2} for the quadratics nonlinear Schr\"odinger equations
\begin{eqnarray}\label{QNLS1}
iu_{t}+u_{xx}+F_j(u,\bar{u})=0, \,\,\, j=1,2,3.
\end{eqnarray}

Here $\bar{u}$ denotes the complex conjugate of $u$ and $F_1(u,\bar{u})=u^2$, $F_2(u,\bar{u})=u\bar{u}$, $F_3(u,\bar{u})={\bar{u}}^2$ in one spatial dimension and in spatially real and periodic case.  We should mention that Bejenaru and Tao \cite{BT} improved the result in \cite{KPV2} for nonlinearity $F_1$ in the real case.

The arguments in \cite {KPV2} use some arithmetic facts involving the symbol of the  linearized equation. For example, the algebraic relation for quadratic NLS (\ref{QNLS1}) with $j=1$ is given by
\begin{equation}\label{AR}
2|n_1(n-n_1)|\leq |\tau-n^2|+|(\tau-\tau_1)-(n-n_1)^2|+|\tau_1-n_1^2|.
\end{equation}

Then splitting the domain of integration in the sets where each term on the right side of (\ref{AR}) is the biggest one, Kenig, Ponce and Vega made some cancellation in the symbol in order to use his calculus inequalities (see Lemma \ref{l3.1}) and a clever change of variables to established their crucial estimates.

To describe our results we define next the $X_{s,b}$ spaces related to our problem.

\begin{defi}\label{GAM}
Let $\mathcal{Y}$ be the space of functions $F(\cdot)$ such that
\begin{enumerate}
\item [$(i)$] $F: \mathbb{T}\times \R \rightarrow \C$.
\item [$(ii)$] $F(x,\cdot)\in S(\R)$ for each $x\in \mathbb{T}$.
\item [$(iii)$] $F(\cdot,t)\in C^{\infty}(\mathbb{T})$ for each $t\in \R$.
\end{enumerate}
For $s,b \in \R$, $X_{s,b}$ denotes the completion of $\mathcal{Y}$ with respect to the norm
\begin{eqnarray}
\|F\|_{X_{s,b}}&=&\|\langle|\tau|-\gamma(n)\rangle^b\langle n\rangle^s \widetilde{F}\|_{l^{2}_{n}L^{2}_{\tau}},
\end{eqnarray}
where $\sim$ denotes the time-space Fourier transform and $\gamma(n)\equiv\sqrt{{n}^2+{n}^4}$.
\end{defi}

Now we state the main results of this paper.

\begin{theo}\label{t1.1}
Let $s >-1/4$ and $u,v\in X_{s,-a}$. Then, there exists $c>0$ depending only on $a,b$ and $s$ such that
\begin{equation}\label{BE}
\left\|
\left(\dfrac{|n|^2\widetilde{uv}(\tau,n)}{2i\gamma(n)}\right)^{\sim^{-1}}
\right\|_{X_{s,-a}}\leq c\left\|u\right\|_{X_{s,b}}\left\|v\right\|_{X_{s,b}},
\end{equation}
where $\sim^{-1}$ denotes the inverse time-space Fourier transform, holds in the following cases
\begin{enumerate}
\item [($i$)] $s\geq 0$, $b>1/2$ and $1/4<a<1/2$,
\item [($ii$)]$-1/4<s<0$, $b>1/2$ and $1/4<a<1/2$ such that $|s|<a/2$.
\end{enumerate}
\end{theo}

\begin{theo}\label{t1.2}
For any $s< -1/4$ and any $a, b\in\R$, with $a<1/2$ the estimate (\ref{BE}) fails. 
\end{theo}

Theorem \ref{t1.2} implies that in our frame the following local well-posed result is sharp (except for the limiting case which remains open).

\begin{theo}\label{t1.3}
Let $s>-1/4$, then for all $\phi\in H^s(\T)$ and $\psi\in H^{s-1}(\T)$, there exist $T=T(\|\phi\|_{H^s},\|\psi\|_{H^{s-1}})$ and a unique solution $u$ of the IVP (\ref{NLB}) with $f(u)=u^2$, $u_0=\phi$ and $u_1=\psi_x$ such that
\begin{equation*}
u\in C([0,T]:H^s(\T))\cap X_{s,b}.
\end{equation*}

Moreover, given $T'\in (0,T)$ there exists $R=R(T')>0$ such that giving the set $ W\equiv\{(\tilde{\phi},\tilde{\psi})\in H^s(\T)\times H^{s-1}(\T):\|\tilde{\phi}-\phi\|_{H^s(\T)}^2+ \|\tilde{\psi}-\psi\|_{H^{s-1}(\T)}^2<R\}$ the map solution 
\begin{equation*}
S:W \longrightarrow C([0,T']:H^s(\T))\cap X_{s,b}, \peq (\tilde{\phi},\tilde{\psi})\longmapsto u(t)
\end{equation*}
is Lipschitz.

In addition, if $(\phi,\psi)\in H^{s'}(\T)\times H^{s'-1}(\T)$ with $s'>s$, then the above results hold with $s'$ instead of $s$ in the same interval $[0,T]$ with 
\begin{equation*}
T=T(\|\phi\|_{H^s},\|\psi\|_{H^{s-1}}).
\end{equation*} 
\end{theo}

These theorems are surprising. The Fourier restriction norm approach, in most of the equations, provide different best Sobolev index results in the real and periodic cases. By means of this technique we have the following LWP results.

For the Korteweg-de Vries equation (KdV)
\begin{eqnarray*}
u_{t}+u_{xxx}+uu_x=0,
\end{eqnarray*}
$s>-3/4$ for the IVP and $s>-1/2$ for the PBVP (see Kenig, Ponce and Vega \cite{KPV1}).

For the modified Korteweg-de Vries equation (mKdV)
\begin{eqnarray*}
u_{t}+u_{xxx}+6u^2u_x=0,
\end{eqnarray*}
$s\geq1/4$ for the IVP and $s\geq1/2$ for the PBVP (see Kenig, Ponce and Vega \cite{KPV6} and Bourgain \cite{B}).

For the equation (\ref{QNLS1}) with nonlinearity $F_1$, $s>-1$ for the IVP and $s>-1/2$ for the PBVP (see Bejenaru and Tao \cite{BT} and Kenig, Ponce and Vega \cite{KPV2}).

For the equation (\ref{QNLS1}) with nonlinearity $F_3$, resp. $F_2$, $s>-3/4$ for the IVP and $s>-1/2$ for the PBVP, resp. $s>-1/4$ for the IVP and $s>0$ for the PBVP (see Kenig, Ponce and Vega \cite{KPV2}).


So it seems that in general the lower Sobolev index which guarantee LWP for the IVP is lower (in fact, at least $1/4$ for the one dimensional case) than that for the IBVP.

In our case the situation is different. Together with Theorem $1.1$-$1.3$ of \cite{LG4}, we conclude that the ``good" Boussinesq equation (\ref{NLB}) in the real and periodic cases are local well-posed in all Sobolev spaces $H^s$ for $s>-1/4$. Since we also show that our estimates for the bilinear operators fail for $s<-1/4$, the result concerning the local well posedness of the ``good" Boussinesq equation is the optimal one provided by the method in both settings (except for the limiting case). As far we know, this is the first example where we have this kind of result for negative indices of $s$.

The plan of this paper is as follows: in Section 2, we prove some estimates for the integral equation in the $X_{s,b}$ space introduced above. Bilinear estimates and the relevant counterexamples are proved in Section 3 and 4, respectively.



\section{Preliminary Results}

First, we consider the linear equation
\begin{equation}\label{LB} 
u_{tt}-u_{xx}+u_{xxxx}=0
\end{equation} 
the solution for initial data $u(0)=\phi$ and $u_t(0)=\psi_x$, is given by
\begin{equation}\label{GUB} 
u(t)=V_c(t)\phi+V_s(t)\psi_x
\end{equation}
where
\begin{eqnarray*}
V_c(t)\phi&=&\left( \frac{e^{it\sqrt{{n}^2+{n}^4}}+ e^{-it\sqrt{{n}^2+{n}^4}}}{2}\hat{\phi}(n)\right)^{\vee},\\
V_s(t){\psi_x}&=&\left( \frac{e^{it\sqrt{{n}^2+{n}^4}}- e^{-it\sqrt{{n}^2+{n}^4}}}{2i\sqrt{{n}^2+{n}^4}}\hat{\psi_x}(n)\right)^{\vee}.
\end{eqnarray*}

By Duhamel's Principle the solution of (NLB) is equivalent to
\begin{equation}\label{INT} 
u(t)= V_c(t)\phi+V_s(t)\psi_x+\int_{0}^{t}V_s(t-t')(u^2)_{xx}(t')dt'.
\end{equation}

Let $\theta$ be a cutoff function satisfying $\theta \in C^{\infty}_{0}(\R)$, $0\leq \theta \leq 1$, $\theta \equiv 1$ in $[-1,1]$, supp$(\theta) \subseteq [-2,2]$ and for $0<T<1$ define $\theta_T(t)=\theta(t/T)$. In fact, to work in the $X_{s,b}$ spaces we consider another version of (\ref{INT}), that is
\begin{equation}\label{INT2} 
u(t)= \theta(t)\left(V_c(t)\phi+V_s(t)\psi_x\right)+\theta_T(t)\int_{0}^{t} V_s(t-t')(u^2)_{xx}(t')dt'.
\end{equation}

Note that the integral equation (\ref{INT2}) is defined for all $(t,x)\in \R^2$. Moreover if $u$ is a solution of (\ref{INT2}) than $\tilde{u}=u|_{[0,T]}$ will be a solution of (\ref{INT}) in $[0,T]$.\\

In the next lemma, we estimate the linear part of the integral equation (\ref{INT2}).

\begin{lemm}\label{l21}
Let $u(t)$ the solution of the linear equation 
\begin{eqnarray*}
\left\{ 
\begin{array}{l}
u_{tt}-u_{xx}+u_{xxxx}=0,\\
u(0,x)=\phi(x); \peq u_t(0,x)=(\psi(x))_x
\end{array} \right. 
\end{eqnarray*}
with $\phi \in H^s$ and $\psi \in H^{s-1}$. Then there exists $c>0$ depending only on $\theta,s,b$ such that
\begin{equation}\label{LP}
\|\theta u\|_{X_{s,b}}\leq c\left(\|\phi\|_{H^s}+\|\psi\|_{H^{s-1}}\right).
\end{equation}
\end{lemm}

\noindent\textbf{Proof. } The proof is analogous to Lemma $2.1$ of \cite{LG4}.\\
\fim

Next we estimate the integral part of (\ref{INT2}).

\begin{lemm}\label{L22}
Let $-\frac{1}{2}<b'\leq 0\leq b \leq b'+1$ and $0<T \leq 1$ then
$$\left\|\theta_T(t)\int_{0}^{t}V_s(t-t')g(u)(t')dt'\right\|_{X_{s,b}} \leq T^{1-(b-b')} \left\|
\left(\dfrac{\widetilde{g(u)}(\tau,n)}{2i\gamma(n)}\right)^{\sim^{-1}}
\right\|_{X_{s,b'}}.$$
\end{lemm}

\noindent\textbf{Proof. } The proof is analogous to Lemma $2.2$ of \cite{LG4}.\\
\fim

To finish this section, we remark that for any positive numbers $a$ and $b$, the notation $a \lesssim b$ means that there exists a positive
constant $\theta$ such that $a \leq \theta b$. We also denote $a \sim b$ when, $a \lesssim b$ and $b \lesssim a$.

\section{Bilinear estimates}

Before proceed to the proof of Theorem \ref{t1.1}, we state some elementary calculus inequalities that will be useful later.
\begin{lemm}\label{l3.1}
For $p, q>0$ and $r=\min\{ p, q, p+q-1\}$ with $p+q>1$, there exists $c>0$ such that
\begin{equation}\label{CI1}
\int_{-\infty}^{+\infty}\dfrac{dx} {\langle x-\alpha\rangle^{p}\langle x-\beta\rangle^{q}}\leq\dfrac{c} {\langle \alpha-\beta\rangle^{r}}.
\end{equation}
\end{lemm}

\noindent\textbf{Proof. } See Lemma 4.2 in \cite{GTV}.\\
\fim

\begin{lemm}\label{l3.2}
If $\gamma>1/2$, then
\begin{equation}\label{CI2}
\sup_{(n,\tau)\in \Z\times\R}\sum_{n_1\in\Z}\dfrac{1}{(1+|\tau\pm n_1(n-n_1)|)^{\gamma}}<\infty. 
\end{equation}
\end{lemm}

\noindent\textbf{Proof. } See Lemma 5.3 in \cite{KPV2}.\\
\fim

\begin{lemm}\label{l3.4}
Let $0<a<1/2$, $\alpha \in \R$, $\beta, \nu>0$ and $H=\{h\in \R: h=\alpha\pm n, n\in \Z \textrm{ and } |h|\leq \beta\}$. Then
\begin{equation}\label{LGI}
\sum_{h\in H}\dfrac{1}{(\nu+|h|)^{2a}}\leq 2\left(\dfrac{2}{\nu^{2a}}+\int_0^{\beta}\dfrac{dx}{(\nu+x)^{2a}} \right). 
\end{equation}
\end{lemm}

\noindent\textbf{Proof. } Since this is only a calculation we omit the proof.\\
\fim

\begin{lemm}\label{l3.3}
There exists $c>0$ such that
\begin{equation}\label{LN}
\dfrac{1}{c}\leq\sup_{x,y\geq 0}\dfrac{1+|x-y|}{1+|x-\sqrt{y^2+y}|}\leq c.
\end{equation}
\end{lemm}

\noindent\textbf{Proof. } Since $y\leq\sqrt{y^2+y}\leq y+1/2$ for all $y\geq 0$ a simple computation shows the desired inequalities.\\
\fim

\begin{rema}\label{R1}
In view of the previous lemma we have an equivalent way to compute the $X_{s,b}$-norm, that is 
\begin{equation*}
\|u\|_{X_{s,b}}\sim\|\langle|\tau|-n^2\rangle^b\langle n\rangle^{s} \widetilde{u}(\tau,n)\|_{l^{2}_{n}L^2_{\tau}}. 
\end{equation*}
This equivalence will be important in the proof of Theorem \ref{t1.1}. As we commented in the introduction the Boussinesq symbol $\sqrt{{n}^2+{n}^4}$ does not have good cancellations to make use of Lemma \ref{l3.1}. Therefore, we  modify the symbols as above and work only with the algebraic relations for the Schr\"odinger equation already used in Kenig, Ponce and Vega \cite{KPV2} in order to derive the bilinear estimates. 
\end{rema}

Now we are in position to prove the bilinear estimate (\ref{BE}).\\

\noindent\textbf{Proof of Theorem \ref{t1.1}.}
Let $u, v\in X_{s,b}$ and define $f(\tau,n)\equiv \langle|\tau|-n^2\rangle^b\langle n\rangle^{s} \widetilde{u}(\tau,n)$, $g(\tau,n)\equiv \langle|\tau|-n^2\rangle^b\langle n\rangle^{s} \widetilde{v}(\tau,n)$. Using Remark \ref{R1} and a duality argument the desired inequality is equivalent to
\begin{equation}\label{DUA}
\left|W(f,g,\phi)\right|\leq c\|f\|_{l^{2}_{n}L^2_{\tau}}\|g\|_{l^{2}_{n}L^2_{\tau}} \|\phi\|_{l^{2}_{n}L^2_{\tau}}
\end{equation}
where
\begin{eqnarray*}
W(f,g,\phi)&=&  \sum_{n,n_1} \int_{\R^2} \dfrac{|n|^2}{\gamma(n)} \dfrac{\langle n\rangle^{s}}{\langle n_1\rangle^{s} \langle n-n_1\rangle^{s} }\\
&&\times \dfrac{g(\tau_1,n_1)f(\tau-\tau_1,n-n_1) \bar{\phi}(\tau,n)}{\langle|\tau|-n^2\rangle^{a} \langle|\tau_1|-n_1^2\rangle^{b} \langle|\tau-\tau_1|-(n-n_1)^2\rangle^{b} } d\tau d\tau_1.
\end{eqnarray*}

Therefore to perform the desired estimate we need to analyze all the possible cases for the sign of $\tau$, $\tau_1$ and $\tau-\tau_1$. To do this we split $\R^4$ into the regions 
\begin{eqnarray*}
\Gamma_1&=&\{(n, \tau, n_1, \tau_1)\in \R^4: \tau_1< 0, \tau-\tau_1< 0\},\\ 
\Gamma_2&=&\{(n, \tau, n_1, \tau_1)\in \R^4: \tau_1\geq 0, \tau-\tau_1< 0, \tau\geq 0\},\\ 
\Gamma_3&=&\{(n, \tau, n_1, \tau_1)\in \R^4: \tau_1\geq 0, \tau-\tau_1< 0, \tau< 0\},\\ 
\Gamma_4&=&\{(n, \tau, n_1, \tau_1)\in \R^4: \tau_1< 0, \tau-\tau_1\geq 0, \tau\geq 0\},\\
\Gamma_5&=&\{(n, \tau, n_1, \tau_1)\in \R^4: \tau_1< 0, \tau-\tau_1\geq 0, \tau< 0\},\\
\Gamma_6&=&\{(n, \tau, n_1, \tau_1)\in \R^4: \tau_1\geq 0, \tau-\tau_1\geq 0\}.
\end{eqnarray*}

Thus, it is sufficient to prove inequality (\ref{DUA}) with $Z(f,g,\phi)$ instead of $W(f,g,\phi)$, where
 \begin{equation*}
Z(f,g,\phi)=  \sum_{n,n_1} \int_{\R^2} \dfrac{|n|^2}{\gamma(n)} \dfrac{\langle n\rangle^{s}}{\langle n_1\rangle^{s} \langle n_2\rangle^{s}} \dfrac{g(\tau_1, n_1)f(\tau_2, n_2) \bar{\phi}(\tau, n)}{\langle\sigma\rangle^{a} \langle\sigma_1\rangle^{b} \langle\sigma_2\rangle^{b}}d\tau d\tau_1
\end{equation*}
with $n_2=n-n_1$, $\tau_2=\tau-\tau_1$ and $\sigma, \sigma_1, \sigma_2$ belonging to one of the following cases
\begin{enumerate}
\item [$(I)$] $\sigma=\tau+n^2,\peq \sigma_1=\tau_1+n_1^2,\peq \sigma_2=\tau_2+n_2^2$,
\item [$(II)$]$\sigma=\tau-n^2,\peq \sigma_1=\tau_1-n_1^2,\peq \sigma_2=\tau_2+n_2^2$,
\item [$(III)$] $\sigma=\tau+n^2,\peq \sigma_1=\tau_1-n_1^2,\peq \sigma_2=\tau_2+n_2^2$,
\item [$(IV)$] $\sigma=\tau-n^2,\peq \sigma_1=\tau_1+n_1^2,\peq \sigma_2=\tau_2-n_2^2$,
\item [$(V)$] $\sigma=\tau+n^2,\peq \sigma_1=\tau_1+n_1^2,\peq \sigma_2=\tau_2-n_2^2$,
\item [$(VI)$] $\sigma=\tau-n^2,\peq \sigma_1=\tau_1-n_1^2,\peq \sigma_2=\tau_2-n_2^2$.
\end{enumerate}

\begin{rema}
Note that the cases $\sigma=\tau+n^2,\peq \sigma_1=\tau_1-n_1^2,\peq \sigma_2=\tau_2-n_2^2$ and $\sigma=\tau-n^2,\peq \sigma_1=\tau_1+n_1^2,\peq \sigma_2=\tau_2+n_2^2$ cannot occur, since $\tau_1< 0, \tau-\tau_1< 0$ implies $\tau<0$ and $\tau_1\geq 0, \tau-\tau_1\geq 0$ implies $\tau\geq 0$.
\end{rema}
Applying the change of variables $(n, \tau, n_1, \tau_1)\mapsto -(n, \tau, n_1, \tau_1)$ and observing that the $l^2_nL^2_{\tau}$-norm is preserved under the reflection operation, the cases $(III)$, $(II)$, $(I)$ can be easily reduced, respectively, to $(IV)$, $(V)$, $(VI)$. Moreover,  making the change of variables $\tau_2=\tau-\tau_1$, $n_2=n-n_1$ and then $(n, \tau, n_2, \tau_2)\mapsto -(n, \tau, n_2, \tau_2)$ the case $(V)$ can be reduced $(IV)$. Therefore we need only establish cases $(IV)$ and $(VI)$.  

We first treat the inequality (\ref{DUA}) with $Z(f,g,\phi)$ in the case ($VI$). We should remark that this estimate is exactly inequality $(1.27)$ that appear in Theorem $1.8$ of \cite{KPV2}, but since it is important to have the inequality (\ref{DUA}) with $a<1/2<b$ such that $a+b<1$ to make the contraction arguments work we reprove this inequalities here.  We will make use of the following algebraic relation
\begin{equation}\label{AR1}
 -(\tau-n^2)+(\tau_1-n_1^2)+((\tau-\tau_1)-(n-n_1)^2)=2n_1(n-n_1).
\end{equation}

By symmetry we can restrict ourselves to the set
\begin{equation*}
 A=\{(n, n_1, \tau, \tau_1)\in \Z^2\times \R^2: |(\tau-\tau_1)-(n-n_1)^2|\leq|\tau_1-n_1^2|\}.
\end{equation*}

We divide $A$ into three pieces
\begin{eqnarray*}
A_1&=&\{(n, n_1, \tau, \tau_1)\in A: n=0\},\\
A_2&=&\{(n, n_1, \tau, \tau_1)\in A: n_1=0 \textrm{ or } n_1 = n \},\\
A_3&=&\{(n, n_1, \tau, \tau_1)\in A: n \neq 0, n_1 \neq 0 \textrm{ and } n_1 \neq n \}.
\end{eqnarray*} 

Next we split $A_3$ into two parts
\begin{eqnarray*}
A_{3,1}&=&\{(n, n_1, \tau, \tau_1)\in A_3: |\tau_1-n_1^2|\leq|\tau-n^2|\},\\
A_{3,2}&=&\{(n, n_1, \tau, \tau_1)\in A_3: |\tau-n^2|\leq|\tau_1-n_1^2| \}.
\end{eqnarray*} 

We can now define the sets $R_i$, $i=1,2$, as follows
\begin{equation*}
R_1=A_1\cup A_2\cup A_{3,1} \textrm{  and  } R_2=A_{3,2}.
\end{equation*}

In what follows $\chi_R$ denotes the characteristic function of the set $R$. Using the Cauchy-Schwarz and H\"older inequalities it is easy to see that
\begin{eqnarray*}
|Z|^2&\leq& \|f\|_{l^{2}_{n}L^2_{\tau}}^2 \|g\|_{l^{2}_{n}L^2_{\tau}}^2 \|\phi\|_{l^{2}_{n}L^2_{\tau}}^2\\
&&\times\left\|\dfrac{1}{\langle\sigma\rangle^{2a}}
\dfrac{|n|^4}{\gamma(n)^2}  \sum_{n_1}\int\dfrac{\langle n\rangle^{2s} \langle n_1\rangle^{-2s} \langle n_2\rangle^{-2s} \chi_{R_1} }{ \langle\sigma_1\rangle^{2b} \langle\sigma_2\rangle^{2b}}d\tau_1\right\|_{l^{\infty}_nL^{\infty}_{\tau}}\\
&&+\|f\|_{l^{2}_{n}L^2_{\tau}}^2 \|g\|_{l^{2}_{n}L^2_{\tau}}^2 \|\phi\|_{l^{2}_{n}L^2_{\tau}}^2\\
&&\times\left\|\dfrac{1}{\langle\sigma_1\rangle^{2b}} \sum_{n}\int \dfrac{|n|^4}{\gamma(n)^2} \dfrac{\langle n\rangle^{2s}\langle n_1\rangle^{-2s} \langle n_2\rangle^{-2s}\chi_{R_2}} {\langle\sigma\rangle^{2a}   \langle\sigma_2\rangle^{2b}}d\tau \right\|_{l^{\infty}_{n_1}L^{\infty}_{\tau_1}}.\\
\end{eqnarray*}

Therefore in view of Lemma \ref{l3.1}-(\ref{CI1}) it suffices to get bounds for
\begin{eqnarray*}
J_1&\equiv&\sup_{n,\tau}\dfrac{1}{\langle\sigma\rangle^{2a}}  \dfrac{|n|^4}{\gamma(n)^2}\sum_{n_1}\dfrac{\langle n\rangle^{2s} \langle n_1\rangle^{-2s} \langle n_2\rangle^{-2s}} {\langle\tau-n^2-2n_1^2+2nn_1\rangle^{2b}}\peq \textrm{on} \peq R_1,\\
J_2&\equiv&\sup_{n_1,\tau_1}\dfrac{1} {\langle\sigma_1\rangle^{2b}} \sum_{n}\dfrac{|n|^4}{\gamma(n)^2}\dfrac{\langle n\rangle^{2s} \langle n_1\rangle^{-2s} \langle n_2\rangle^{-2s}} {\langle\tau_1+n_1^2-2nn_1\rangle^{2a}}\peq \textrm{on} \peq R_2.
\end{eqnarray*}

In the region $A_1$ we have $\dfrac{|n|^4}{\gamma(n)^2}=0$, therefore the estimate is trivial.\\

In region $A_2$, we have 
\begin{equation*}
\langle n\rangle^{2s} \langle n_1\rangle^{-2s} \langle n_2\rangle^{-2s}\lesssim 1 \peq \textrm{and} \peq \dfrac{|n|^4}{\gamma(n)^2}\leq 1 \peq \textrm{for all} \peq n\in\Z.
\end{equation*}

Therefore for $n_1=0$ or $n_1=n$ we obtain
\begin{equation*}
J_1 \lesssim \sup_{n,\tau} \dfrac{1}{\langle\sigma\rangle^{2a+2b}} \lesssim 1
\end{equation*}
for $a,b>0$.\\

Now, by definition of region $A_{3,1}$ and the algebraic relation (\ref{AR1}) we have
\begin{equation*}
\langle n\rangle^{2s} \langle n_1\rangle^{-2s} \langle n_2\rangle^{-2s} \lesssim \langle \sigma\rangle^{\lambda(s)}
\end{equation*}
where 
\begin{eqnarray}\label{LAMBDA}
\lambda(s)=
\left\{ 
\begin{array}{l c}
0, &\textrm{ if } s\geq 0\\
2|s|, &\textrm{ if } s\leq 0
\end{array} \right. .
\end{eqnarray}

Therefore by Lemma \ref{l3.2}-(\ref{CI2})
\begin{equation*}
J_1 \lesssim \sup_{n,\tau} \langle\sigma\rangle^{\lambda(s)-2a}\sum_{n_1\neq \{0,n\}} \dfrac{1}{\langle\tau-n^2-2n_1^2+2nn_1\rangle^{2b}} 
\lesssim 1
\end{equation*}
for $b>1/2$ and ${\lambda(s)}\leq 2a$.\\

Next we estimate $J_2$. In region $A_{3,2}$, we have
\begin{equation*}
\langle n\rangle^{2s} \langle n_1\rangle^{-2s} \langle n_2\rangle^{-2s} \lesssim \langle \sigma_1\rangle^{\lambda(s)}
\end{equation*}
and
\begin{equation*}
|\tau_1+{n_1}^2-2nn_1|\leq 2\langle\sigma_1\rangle.
\end{equation*} 

Define $H=\{n\in \Z: |\tau_1+{n_1}^2-2nn_1|\leq 2\langle\sigma_1\rangle\}$. In view of Lemma \ref{l3.4}-\eqref{LGI}, we have 
\begin{eqnarray*}
J_2 &\lesssim& \sup_{n_1,\tau_1}
\langle \sigma_1\rangle^{\lambda(s)-2b}\sum_{n\in H} \dfrac{1}{\langle\tau_1+{n_1}^2-2nn_1\rangle^{2a}}\\
&\lesssim& \sup_{n_1,\tau_1} \dfrac{\langle \sigma_1\rangle^{\lambda(s)-2b} }{|n_1|^{2a}}\left(|n_1|^{2a}+ \int_{0 }^{\langle\sigma_1\rangle/|n_1|} {\left(\dfrac{1}{2|n_1|} +x\right)^{-2a}}dx \right)\\
&\lesssim& \sup_{n_1,\tau_1} \langle \sigma_1\rangle^{\lambda(s)-2b} + \dfrac{\langle \sigma_1\rangle^{\lambda(s)-2b} }{|n_1|} + \dfrac{\langle \sigma_1\rangle^{\lambda(s)-2a -2b+1} }{|n_1|}
\lesssim 1
\end{eqnarray*}
for $a<1/2$ and $\lambda(s)\leq \min\{2b,2a+2b-1\}$.\\

Now we turn to the proof of case ($IV$). This is analogous to inequality $(1.29)$ of \cite{KPV2}. In the following estimates we will make use of the algebraic relation
\begin{equation}\label{AR2}
 -(\tau-n^2)+(\tau_1+n_1^2)+((\tau-\tau_1)-(n-n_1)^2)=2n_1n.
\end{equation}

First we split $\Z^2\times \R^2$ into three sets
\begin{eqnarray*}
B_1&=&\{(n, n_1, \tau, \tau_1)\in \Z^2\times \R^2: n=0\},\\
B_2&=&\{(n, n_1, \tau, \tau_1)\in \Z^2\times \R^2: n_1=0\},\\
B_3&=&\{(n, n_1, \tau, \tau_1)\in \Z^2\times \R^2: n_1\neq 0, n\neq 0\}.
\end{eqnarray*} 

Next we separate $B_3$ into three parts
\begin{small}
\begin{eqnarray*}
B_{3,1}&=&\{(n, n_1, \tau, \tau_1)\in B_3: |\tau_1+n_1^2|,|(\tau-\tau_1)-(n-n_1)^2|\leq|\tau-n^2|\},\\
B_{3,2}&=&\{(n, n_1, \tau, \tau_1)\in B_3: |\tau-n^2|,|(\tau-\tau_1)-(n-n_1)^2|\leq|\tau_1+n_1^2|\},\\
B_{3,3}&=&\{(n, n_1, \tau, \tau_1)\in B_3: |\tau_1+n_1^2|,|\tau-n^2|\leq|(\tau-\tau_1)-(n-n_1)^2|\}.
\end{eqnarray*} 
\end{small}

We can now define the sets $S_i$, $i=1,2$, as follows
\begin{equation*}
S_1=B_1\cup B_2\cup B_{3,1}, \peq S_2=B_{3,2} \peq \textrm{and} \peq S_3=B_{3,3}.
\end{equation*}

Using the Cauchy-Schwarz and H\"older inequalities and duality it is easy to see that
\begin{eqnarray*}
|Z|^2\!\!\!&\lesssim&\!\!\! \|f\|_{l^{2}_{n}L^2_{\tau}}^2 \|g\|_{l^{2}_{n}L^2_{\tau}}^2 \|\phi\|_{l^{2}_{n}L^2_{\tau}}^2\\
&&\times\left\|\dfrac{1}{\langle\sigma\rangle^{2a}} 
\dfrac{|n|^4}{\gamma(n)^2}
\sum_{n_1}\int\dfrac{\langle n\rangle^{2s} \langle n_1\rangle^{-2s} \langle n_2\rangle^{-2s} \chi_{S_1} }{\langle\sigma_1\rangle^{2b} \langle\sigma_2\rangle^{2b}}d\tau_1\right\|_{l^{\infty}_nL^{\infty}_{\tau}}\\
&&+\|f\|_{l^{2}_{n}L^2_{\tau}}^2 \|g\|_{l^{2}_{n}L^2_{\tau}}^2 \|\phi\|_{l^{2}_{n}L^2_{\tau}}^2\\
&&\times\left\|\dfrac{1}{\langle\sigma_1\rangle^{2b}} \sum_{n}\int \dfrac{|n|^4}{\gamma(n)^2} \dfrac{\langle n\rangle^{2s} \langle n_1\rangle^{-2s} \langle n_2\rangle^{-2s} \chi_{S_2}}{\langle\sigma\rangle^{2a}   \langle\sigma_2\rangle^{2b}}d\tau\right\|_{l^{\infty}_{n_1}L^{\infty}_{\tau_1}}\\
&&+\|f\|_{l^{2}_{n}L^2_{\tau}}^2 \|g\|_{l^{2}_{n}L^2_{\tau}}^2 \|\phi\|_{l^{2}_{n}L^2_{\tau}}^2\\
&&\times\left\|\dfrac{1}{\langle\sigma_2\rangle^{2b}} \sum_{n_1}\int \dfrac{|n_1+n_2|^4}{\gamma(n_1+n_2)^2} \dfrac{\langle n_1+ n_2\rangle^{2s} \langle n_1\rangle^{-2s} \langle n_2\rangle^{-2s} \chi_{\widetilde{S}_3}}{\langle\sigma_1\rangle^{2b}   \langle\sigma\rangle^{2a}}d\tau_1\right\|_{l^{\infty}_{n_2}L^{\infty}_{\tau_2}}\\
\end{eqnarray*}
where $\sigma$, $\sigma_1$, $\sigma_2$ were given in the condition ($IV$) and
\begin{eqnarray*}
\widetilde{S}_3\subseteq 
\left\{ 
\begin{array}{r}
(n_2, n_1, \tau_2, \tau_1)\in \Z^2\times\R^2: n_1\neq 0, n_1+n_2\neq 0 \textrm{ and }\\
|\tau_1+n_1^2|,|(\tau_1+\tau_2)-(n_1+n_2)^2|\leq|\tau_2-n_2^2|
\end{array} \right\}.
\end{eqnarray*} 

Therefore from Lemma \ref{l3.1}-(\ref{CI1}) it suffices to get bounds for
\begin{eqnarray*}
K_1\!\!\!&\equiv&\!\!\!\sup_{n,\tau} \dfrac{1}{\langle\sigma\rangle^{2a}} 
\dfrac{|n|^4}{\gamma(n)^2} \sum_{n_1}\dfrac{\langle n\rangle^{2s}\langle n_1\rangle^{-2s} \langle n_2\rangle^{-2s}} {\langle\tau-n^2+2nn_1\rangle^{2b}}\peq \textrm{on} \peq S_1,\\
K_2\!\!\!&\equiv&\!\!\!\sup_{n_1,\tau_1}\dfrac{1} {\langle\sigma_1\rangle^{2b}} \sum_{n} \dfrac{|n|^4}{\gamma(n)^2}\dfrac{\langle n\rangle^{2s} \langle n_1\rangle^{-2s} \langle n_2\rangle^{-2s}} {\langle\tau_1+n_1^2-2nn_1\rangle^{2a}}\peq \textrm{on} \peq S_2,\\
K_3\!\!\!&\equiv&\!\!\!\sup_{n_2,\tau_2}\dfrac{1}{\langle\sigma_2\rangle^{2b}} \sum_{n_1}\dfrac{|n_1+n_2|^4}{\gamma(n_1+n_2)^2}\dfrac{\langle n_1+n_2\rangle^{2s} \langle n_1\rangle^{-2s} \langle n_2\rangle^{-2s}} {\langle\tau_2-n_2^2-2n_1^2-2n_1n_2\rangle^{2a}}\peq \textrm{on} \peq \widetilde{S}_3.
\end{eqnarray*}

In the region $B_1$ we have $\dfrac{|n|^4}{\gamma(n)^2}=0$, therefore the estimate is trivial.\\

\begin{rema}\label{re0}
We should notice that region $B_1$ is difficult to handle when we consider equation (\ref{QNLS1}) with nonlinearity $F_2$ in the periodic setting (see Kenig, Ponce and Vega \cite{KPV2}, Theorem $1.10-(iii)$). In order to prove the bilinear estimate in this case one need to bound
\begin{equation*}
\widetilde{K}_2\equiv \sup_{n_1,\tau_1}\dfrac{1} {\langle \tau_1+n_1^2 \rangle^{2b}} \sum_{n} \dfrac{\langle n\rangle^{2s} \langle n_1\rangle^{-2s} \langle n_2\rangle^{-2s}} {\langle\tau_1+n_1^2-2nn_1\rangle^{2a}}.
\end{equation*}

In particular, when $n=0$ we have
\begin{equation*}
\widetilde{K}_2\geq \sup_{n_1,\tau_1}\dfrac{\langle n_1\rangle^{-4s}} {\langle\tau_1+n_1^2\rangle^{2b+2a}}.
\end{equation*}

Now, taking $\tau_1=-N^2$, $n_1=N$ and letting $N\rightarrow \infty$ we conclude $\widetilde{K}_2=\infty$ for $s<0$. Therefore, we cannot push the LWP to equation (\ref{QNLS1}) with nonlinearity $F_2$ for any $s<0$ using this method.
\end{rema}

In region $B_2$, we have 
\begin{equation*}
\langle n\rangle^{2s} \langle n_1\rangle^{-2s} \langle n_2\rangle^{-2s}\lesssim 1 \peq \textrm{and} \peq \dfrac{|n|^4}{\gamma(n)^2}\leq 1 \peq \textrm{for all} \peq n\in\Z.
\end{equation*}

Therefore for $n_1=0$ we obtain
\begin{equation*}
K_1\lesssim \sup_{n,\tau}\dfrac{1}{\langle\sigma\rangle^{2a+2b}} \lesssim 1
\end{equation*}
for $a,b>0$.\\

In region $B_{3,1}$, we have
\begin{equation}\label{EN}
\langle n\rangle^{2s} \langle n_1\rangle^{-2s} \langle n_2\rangle^{-2s} \lesssim |n_1|^{\eta(s)} \lesssim |n_1n|^{\eta(s)}
\end{equation}
where 
\begin{eqnarray*}
\eta(s)=
\left\{ 
\begin{array}{l c}
0, &\textrm{ if } s\geq 0\\
4|s|, &\textrm{ if } s\leq 0
\end{array} \right. .
\end{eqnarray*}

Moreover, by the algebraic relation (\ref{AR2}) we have
\begin{equation*}
| n_1 n|\lesssim \langle\sigma\rangle.
\end{equation*}

Therefore
\begin{equation*}
K_1 \lesssim \sup_{n,\tau} \langle\sigma\rangle^{\eta(s)-2a}\sum_{n_1\neq 0} \dfrac{1}{\langle\tau-n^2+2nn_1\rangle^{2b}} 
\lesssim 1
\end{equation*}
for $b>1/2$ and ${\eta(s)}\leq 2a$.\\

Next we estimate $K_2$. In the region $B_{3,2}$, we have
\begin{equation*}
| n_1 n|\lesssim \langle\sigma_1\rangle
\end{equation*} 
and
\begin{equation*}
|\tau_1+{n_1}^2-2nn_1|\leq 2\langle\sigma_1\rangle.
\end{equation*}  

In view of (\ref{EN}) and Lemma \ref{l3.4}-\eqref{LGI}, we obtain  
\begin{eqnarray*}
K_2
&\lesssim& \sup_{n_1,\tau_1}
\langle \sigma_1\rangle^{\eta(s)-2b}\sum_{n\in H} \dfrac{1}{\langle\tau_1+{n_1}^2-2nn_1\rangle^{2a}}\\
&\lesssim& \sup_{n_1,\tau_1} \dfrac{\langle \sigma_1\rangle^{\eta(s)-2b} }{|n_1|^{2a}}\left(|n_1|^{2a}+ \int_{0 }^{\langle\sigma_1\rangle/|n_1|} {\left(\dfrac{1}{2|n_1|} +x\right)^{-2a}}dx \right)\\
&\lesssim& \sup_{n_1,\tau_1} \langle \sigma_1\rangle^{\eta(s)-2b} + \dfrac{\langle \sigma_1\rangle^{\eta(s)-2b} }{|n_1|} + \dfrac{\langle \sigma_1\rangle^{\eta(s)-2a -2b+1} }{|n_1|}
\lesssim 1
\end{eqnarray*}
for $a<1/2$ and $\eta(s)\leq \min\{2b,2a+2b-1\}$.\\

Finally, we estimate $K_3(n_2,\tau_2)$. In the region $B_{3,3}$, we have 
\begin{equation*}
\langle n_1+n_2\rangle^{2s} \langle n_1\rangle^{-2s} \langle n_2\rangle^{-2s}\lesssim |n_2|^{\eta(s)}.
\end{equation*}

Moreover, the algebraic relation (\ref{AR2}) implies 
\begin{equation*}
|n_2| \lesssim|n_1+n_2|+|n_1|\lesssim |n_1(n_1+n_2)|\lesssim \langle\sigma_2\rangle.
\end{equation*} 

Therefore, Lemma \ref{l3.2}-\eqref{CI2} implies that
\begin{equation*}
K_3\lesssim \sup_{n_2,\tau_2} \langle\sigma_2\rangle^{\eta(s)-2b}\sum_{n_1\neq 0} \dfrac{1}{\langle\tau_2-n_2^2-2n_1^2-2n_1n_2\rangle^{2a}} \lesssim  1
\end{equation*}
for $a>1/4$, $b>1/2$ and $s>-1/4$ which implies $\eta(s)\leq 2b$.\\
\fim
\begin{rema}
Once the bilinear estimates in Theorem \ref{t1.1} have been established, it is a standard matter to conclude the LWP statement of Theorem \ref{t1.3}. We refer the reader to the works \cite{KPV2} and \cite{LG4} for further details.
\end{rema}

\section{Counterexample to the bilinear estimate (\ref{BE})}

\noindent\textbf{Proof of Theorem \ref{t1.2}.}
For $u\in X_{s,b}$ and $v\in X_{s,b}$ we define $f(\tau,n)\equiv \langle|\tau|-\gamma(n)\rangle^b\langle n\rangle^{s} \widetilde{u}(\tau,n)$ and $g(\tau,n)\equiv \langle|\tau|-\gamma(n)\rangle^b\langle n\rangle^{s} \widetilde{v}(\tau,n)$. By Lemma \ref{l3.3} the inequality (\ref{BE}) is equivalent to
\begin{eqnarray}\label{DE}
\left\|\dfrac{|n|^2}{\gamma(n)}\dfrac{\langle n\rangle^{s}} {\langle\sigma\rangle^{a}}
\sum_{n_1}\int 
\dfrac{ f(\tau_1,n_1) g(\tau_2,n_2) d\tau_1}{\langle n_1\rangle^{s}\langle n_2\rangle^{s}\langle\sigma_1\rangle^b \langle\sigma_2\rangle^b}\right\|_{l^{2}_{n}L^2_{\tau}}\lesssim \|f\|_{l^{2}_{n}L^2_{\tau}}\|g\|_{l^{2}_{n}L^2_{\tau}},
\end{eqnarray}
where
\begin{equation*}
n_2=n-n_1, \peq \tau_2=\tau-\tau_1,
\end{equation*}
\centerline{$
\sigma=|\tau|-n^2,\peq \sigma_1=|\tau_1|-n_1^2,\peq \sigma_2=|\tau_2|-n_2^2.
$}

For $N\in\Z$ define

\centerline{$f_N(\tau,n)=a_n\chi((\tau-n^2)/2)$, with $a_n=\left\{ 
\begin{array}{l}
1, \peq n=N, \\
0, \peq \textrm{elsewhere.}
\end{array} \right.$}
and

\centerline{$g_N(\tau,n)=b_n\chi((\tau+n^2)/2)$, with $b_n=\left\{ 
\begin{array}{l}
1, \peq n=1-N, \\
0, \peq \textrm{elsewhere.}
\end{array} \right.$}
where $\chi(\cdot)$ denotes the characteristic function of the interval $[-1,1]$. Thus
$$a_{n_1}b_{n-n_1}\neq 0 \peq \textrm{if and only if} \peq n_1=N \peq \textrm{and} \peq n=1$$
and consequently for $N$ large
\begin{eqnarray*}
\int\!\!\chi((\tau_1-n_1^2)/2)\chi((\tau-\tau_1+(n-n_1)^2)/2) &\!\!\!\!\gtrsim\!\!\!\!& \chi((\tau+(n-n_1)^2-n_1^2))\\
&\!\!\!\!\gtrsim\!\!\!\!& \chi((\tau+1-2N)).
\end{eqnarray*}

Therefore, using the fact that $||\tau|-n^2|\leq\min\{|\tau-n^2|, |\tau+n^2|\}$, inequality (\ref{DE}) implies 
\begin{eqnarray*}
1 \peq \gtrsim \peq \left\|\dfrac{N^{-2s}}{N^a}\chi((\tau+1-2N))\right\|_{L^2_{\tau}}\peq \gtrsim \peq N^{-2s-a}.
\end{eqnarray*}

Letting $N\rightarrow\infty$, this inequality is possible only when $s\geq -a/2$ which yields the result since $a<1/2$.\\
\fim

\bibliographystyle{amsplain}

E-mail: farah@impa.br and scialom@ime.unicamp.br
\end{document}